\documentclass[12pt]{article}
\usepackage[latin2]{inputenc}
\usepackage{amsfonts, marvosym}
\usepackage{amsmath}
\usepackage{enumerate}
\usepackage{amssymb}
\usepackage{theorem}
\usepackage{authblk}
\usepackage{hyperref}
\usepackage{mathtools}

\newcommand{\R}{{\ensuremath{\mathbb{R}}}}
\newcommand{\N}{{\ensuremath{\mathbb{N}}}}

\renewcommand{\P}{\ensuremath{\mathbb{P}}}
\renewcommand{\dj}{d\kern-0.4em\char"16\kern-0.1em}
\newcommand{\E}{\ensuremath{\mathcal{E}}}
\newcommand{\B}{\ensuremath{\mathcal{B}}}

\newcommand{\proof}{\noindent\textbf{Proof.}\ }
\newcommand{\proofof}{\noindent\textbf{Proof of}\ }

\newcommand{\qed}{\hfill\ensuremath{\Box}\\}

\newcommand{\cp}{\ensuremath{\mathrm{Cap}}}

\newcommand{\Hj}{\hyperlink{(H1)}{\bf(H1)}}
\newcommand{\Hd}{\hyperlink{(H2)}{\bf(H2)}}

\newcommand{\EE}{\mathbb{E}}
\newcommand{\F}{\mathcal{F}}
\newcommand{\C}{\mathcal{C}}

\newtheorem{Thm}{Theorem}[section]

\newtheorem{Cor}[Thm]{Corollary}
\newtheorem{Lem}[Thm]{Lemma}

\newtheorem{Def}[Thm]{Definition}

\newtheorem{Rem}[Thm]{Remark}
{\theorembodyfont{\upshape} }
\numberwithin{equation}{section}

\title{Censored symmetric L\' evy processes\footnote{This work was supported by the Croatian Science Foundation under the project 3526.}}
\author{V.~Wagner\footnote{Department of Mathematics, University of Zagreb, 10000 Zagreb, Croatia\\\Letter\quad\href{mailto:wagner@math.hr}{wagner@math.hr}}}

\begin{document}

\maketitle

\begin{abstract}
We examine three equivalent constructions of a censored symmetric purely discontinuous L\' evy process on an open set $D$; via the corresponding Dirichlet form, through the Feynman-Kac transform of the L\' evy process killed outside of $D$ and from the same killed process by the Ikeda-Nagasawa-Watanabe piecing together procedure. By applying the trace theorem on $n$-sets for Besov-type spaces of generalized smoothness associated with complete Bernstein functions satisfying certain scaling conditions, we analyze the boundary behaviour of the corresponding censored L\'evy process and determine conditions under which the process approaches the boundary $\partial D$ in finite time. Furthermore, we prove a stronger version of the 3G inequality and its generalized version for Green functions of purely discontinuous L\'evy processes on $\kappa$-fat open sets. Using this result, we obtain the scale invariant Harnack inequality for the corresponding censored process.
\end{abstract}


\section{Introduction}

Censored L\' evy process on an open set $D$ in $\R^n$ is a process obtained by restricting (censoring) the jumping measure of a purely discontinuous symmetric L\' evy process to $D\times D$, i.e. by suppressing its jumps outside of $D$. Censored stable processes, obtained from the symmetric $\alpha$-stable L\' evy process, have been introduced by Bogdan, Burdzy and Chen in \cite{cen}, where they analyzed their boundary behaviour, as well as several potential-theoretic properties. Censored stable and stable-like processes have been the center of study of several following papers, for example \cite{green_stable}, \cite{hardy}, \cite{fatou}, \cite{bubbles}. The main goal of this paper is to extend some results obtained for censored stable processes in \cite{cen} to a wider class of discontinuous symmetric L\' evy processes, specifically to analyze their boundary behaviour and prove the Harnack inequality. Additionally, we obtain a more general version of the 3G inequality for this class of L\' evy processes, which, together with the generalized 3G inequality, may be of independent interest.

Let $X=(X_t)_{t\geq 0}$ be a purely discontinuous symmetric L\' evy process in $\R^n$ with characteristic exponent $\Psi_X$ and $D\subset \R^n$ an open set. Following the approach in \cite{cen} which deals with the stable process, in Section \ref{section:construction} we define the censored process $Y=(Y_t)_{t\geq 0}$ on $D$ with lifetime $\zeta$ related to $X$ through its associated regular Dirichlet form $(\E,\F)$ on $L^2(D)$ and present two equivalent constructions - through the Feynman-Kac transform of the killed process $X^D=(X_t^D)_{t\geq 0}$ and by the Ikeda-Nagasawa-Watanabe piecing together procedure applied to $X^D$. These three construction methods provide a wide range of analysis techniques which we use throughout the paper. 

From this point on we restrict ourselves to the case when the L\' evy density $J_X$ of $X$ is comparable to the L\' evy density of an isotropic unimodal L\' evy process. Let $j:(0,\infty)\to(0,\infty)$ be a non-increasing function satisfying
\begin{equation}\label{eq:j1}
j(r) \leq c_1 j(r + 1),\  
\end{equation}
for all $r \geq 1$ and some constant $c_1>0$, such that 
\begin{equation}\label{eq:j2}
 \gamma_1^{-1}j(|x|)\leq J_X(x)\leq \gamma_1j(|x|),
\end{equation}
for some $\gamma_1\geq 1$. Such a function $j$ is a radial L\'evy density of an isotropic unimodal L\' evy process with characteristic exponent
\begin{equation}\label{eq:j3}
 \psi(|\xi|)=\int_{\R^n\setminus\{0\}}\left(1-\cos{(x\cdot\xi)}\right)j(|x|)dx,\ \xi\in\R^n. 
\end{equation}  
Note that \eqref{eq:j2} implies that $\Psi_X\asymp\psi(|\cdot|)$. Throughout the paper we will assume that $\psi$ satisfies one or both of the following scaling conditions, \\
{\hypertarget{(H1)}{\bf(H1)}}: There exist constants $0 < \delta_1 \leq \delta_2 < 1$ and $a_1,a_2 > 0$ such that
$$a_1 \lambda^{2\delta_1} \psi(t) \leq \psi(\lambda t) \leq a_2 \lambda^{2\delta_2} \psi(t),\quad \lambda\geq 1,\, t \geq1;$$

{\hypertarget{(H2)}{\bf(H2)}}: There exist constants $0 < \delta_3 \leq \delta_4 < 1$ and $a_3,a_4 > 0$ such that
$$a_3 \lambda^{2\delta_3} \psi(t) \leq \psi(\lambda t) \leq a_4 \lambda^{2\delta_4} \psi(t),\quad \lambda\geq 1,\, t < 1.$$
Under condition \Hj, by \cite[(2.1), (2.2)]{martin} (see also \cite{BOGDAN20143543}), there exists a complete Bernstein function $\phi$ and a constant $\gamma_2\geq1$ such that
\begin{equation}\label{eq:phipsi}
 \gamma_2^{-1}\phi(|\xi|^2)\leq \psi(|\xi|)\leq \gamma_2\phi(|\xi|^2), \ \xi\in\R^n, 
\end{equation}
and the radial L\'evy density $j$ enjoys the following property: for every $R>0$
\begin{equation}\label{eq:j4}
 j(r)\asymp\frac{\phi(r^{-2})}{r^n},\ r\in(0,R).
\end{equation}
Furthermore, by \cite[Lemma 2.1]{martin} every Bernstein function $\phi$ satisfies the following useful inequality,
\begin{equation}
  1\wedge\lambda\leq\frac{\phi(\lambda r)}{\phi(r)} \leq 1 \vee\lambda,\quad \lambda, r>0.\label{eq:phi}
\end{equation}
In order to analyze the behaviour of the censored process $Y$ near the boundary $\partial D$, we consider the reflected process $Y^*$ corresponding to $Y$ introduced in \cite{silverstein} and relate these two processes through their corresponding Dirichlet forms $(\E,\F)$ and $({\mathcal E}^{\text{ref}},\F_\text{a}^{\text{ref}})$. When $D$ is an open $n$-set (see Definition \ref{def:dset}), $({\mathcal E}^{\text{ref}},\F_\text{a}^{\text{ref}})$ is a regular Dirichlet form on $\overline{D}$ and we can interpret the censored process $Y$ as the reflected process $Y^*$ killed upon hitting the boundary $\partial D$. This is shown by associating the Dirichlet forms $(\C,\F^{\R^n})$ and $({\mathcal E}^{\text{ref}},\F_\text{a}^{\text{ref}})$ with the Besov-type space of generalized smoothness $H^{\psi,1}(\R^n)$ and the corresponding trace space on $n$-set $\overline{D}$, respectively. It follows that the question of $Y$ approaching the boundary in finite time is equivalent to the question of ${\mathcal E}^{\text{ref}}$-polarity of the boundary, and can therefore be partially answered in terms of the Hausdorff dimension of the boundary $\partial D$. Denote by $\mathcal H_h$ the Hausdorff $h$-measure and by $\mathcal H^d$ the $d$-dimensional Hausdorff measure. Finally, we give a characterization of ${\mathcal E}^{\text{ref}}$-polar sets in terms of polarity for $X$ and arrive to the main results of this section, Corollary \ref{cor:boundary} and Theorem \ref{thm:boundary}. 
\begin{Thm}\label{thm:boundary}
Suppose that $D\subset\R^n$ is an open $n$-set and $X$ is a purely discontinuous symmetric L\'evy process such that \eqref{eq:j1}, \eqref{eq:j2}, \Hj and \Hd hold. Let $Y$ be the censored process on $D$ related to $X$ with lifetime $\zeta$ and $Y^*$ the corresponding reflected process on $\overline{D}$.
\begin{enumerate}[(i)]
 \item Suppose that $\delta_2 \leq \frac n 2$ and that $\mathcal H_h (\partial D \cap K_m ) < \infty$ for an increasing sequence of Borel sets $K_m$ such that $\cup_{m\in\N}K_m \supset \partial D$, where $h(r) = r
^{n-2\delta_2}$ if $\delta_2 < \frac n 2$ and $h(r) = \max\{\log \frac 2 r , 0\}$ when $\delta_2 = \frac n 2=\frac 1 2$. Then $\cp_X(\partial D)=0$ and therefore $Y = Y^*$ and $Y$ is conservative.

\item If $\mathcal H^d(\partial D)>0$ for some $d>n-2\delta_1\geq 0$ then $\cp_X(\partial D)>0$, Y is a proper subprocess of $Y^*$, $Y$ is transient and
\[
\P_x(Y_{\zeta-}\in\partial D,\zeta<\infty)>0,\quad\forall x\in D.
\]
If $D$ additionally has finite Lebesgue measure, then $Y$ approaches the boundary in finite time almost surely.
\item Suppose $n=1$. If $\delta_3\geq\frac 1 2$ then $Y^*$ is recurrent. If additionally $\delta_1>\frac 1 2$ then $Y$ is transient and
\[
\P_x(Y_{\zeta-}\in\partial D,\zeta<\infty)=1,\quad\forall x\in D.
\]
\end{enumerate}
\end{Thm}

In Section \ref{section:3G} we prove a stronger version of the so called 3G inequality for purely discontinuous symmetric L\' evy processes on bounded $\kappa$-fat open sets, as well as the generalized 3G inequality. The 3G inequality and generalized 3G inequality are essential tools in obtaining sharp two-sided Green function estimates for local and non-local perturbations of symmetric purely discontinuous L\'evy processes, see \cite{sharp} and \cite{twosided}. The goal is to show that for every $R>0$ and every ball $B\subset D$ of radius $r\leq R$ the Green functions of the killed processes $Y^B$ and $X^B$ are comparable, with constants depending only on $R$ and $\Psi_X$. Since the Green function $G_B$ of the killed L\'evy process $X^B$ lacks the exact scaling property exhibited by the $\alpha$-stable process, the following stronger version of the 3G inequality is needed. For notational convenience, define $\Phi(\lambda)=\frac 1{\phi(\lambda^{-2})}$, $\lambda>0$.  
\begin{Thm}{(3G theorem)}\label{harmonic:3GThm}\\
Let $X$ be a purely discontinuous symmetric L\'evy process such that \eqref{eq:j1}, \eqref{eq:j2} and \Hj hold. Let $r>0$, $a>0$ and $\kappa\in\left(0,\frac 1 2 \right]$. There exists a constant $c_2=c_2(r,a,\kappa,\phi,\gamma_1,\gamma_2)>0$ such that 
\begin{equation}
\frac{G_B(x,y)G_B(y,z)}{G_B(x,z)}\leq c_1 \frac{\Phi(|x-y|)\Phi(|y-z|)}{\Phi(|x-z|)}\frac{|x-z|^n}{|x-y|^n|y-z|^n},\ x,y,z\in B, \label{harmonic:3G}
\end{equation}
for every bounded $\kappa$-fat open set $B$ with characteristics $(R,\kappa)$, diam$(B)\leq r$ and $\frac{R}{\text{diam}(B)}\geq a$.
\end{Thm}

Using the results from the previous section and the representation of the censored process as a Feynman-Kac transform of the killed process $X^D$, in Section \ref{section:harnack} we prove the scale invariant Harnack inequality for non-negative harmonic functions of the censored L\'evy process $Y$.

For easier notation, denote by $d$ the diagonal in $\R^n\times\R^n$. For a bounded set $B$ in $\R^n$ let $\text{diam}(B):=\sup\{|x-y|:\,x,y\in B\}$, $d(x,B):=\inf\{|x-y|:\,y\in B\}$ and $\delta_B(x)=d(x,B^c)$, $x\in \R^n$. We say that functions $f$ and $g$ are comparable and denote $f\asymp g$ if there exists a constant $c>1$ such that $c^{-1}g(x)\leq f(x)\leq c g(x)$ for all $x$. 

\section{Construction and boundary behaviour}\label{section:construction}

Let $(\Omega,\mathcal G,\mathbb P)$ be a probability space and $X=(X_t)_{t\geq0}$ be a purely discontinuous symmetric L\' evy process in $\R^n$ with the L\' evy density $J_X$. The Fourier transform of the transition probability of $X$ is characterized by the L\' evy-Khintchine exponent $\Psi_X(\xi)=\int_{\R^n\setminus\{0\}}\left(1-\cos{(x\cdot\xi)}\right)J_X(x)dx$, $\xi\in\R^n$, 
\begin{align*}
&\mathbb E\left[e^{i\xi\cdot X_t}\right]=
e^{-t\Psi_X(\xi)}.
\end{align*}
The regular Dirichlet form $(\C,\F^{\R^n})$ associated with $X$ is given by
\begin{align}
&\C(u,v)=\frac
1 2\int_{\R^n}\int_{\R^n\setminus\{0\}}(u(x+y)-u(x))(v(x+y)-v(x))J_X(y)dy\,dx\nonumber\\
&\F^{\R^n}=\left\{u\in L^2(\R^n): \C(u,u)<\infty\right\},\nonumber
\end{align}
with $C_c^{\infty}(\R^n)$ as a special standard core, see \cite[Example 1.4.1]{fukushima}. Let $D\subset\R^n$ be an open set and $\tau_D=\inf\{t>0:X_t\not\in D\}$ be the first exit time of $X$ from $D$. Let $X^D=(X^D_t)_{t\geq 0}$ be the process $X$ killed upon exiting $D$, that is 
\[
X_t^D=\left\{\begin{array}{ll}
        X_t, & t\le \tau_D\\
        \partial, & t\geq \tau_D
        \end{array}\right., 
\]
where $\partial$ is the so-called cemetery state. The associated Dirichlet form for $X^D$ is $(\C,\F^D)$, where $\F^D=\{u\in\F^{\R^n}:u=0\ \C\text{-q.e. on }D^c\}$. Here a statement holds $\C$-\emph{quasi-everywhere} ($\C$-\emph{q.e.}) if it holds outside of some set of $\C$-capacity zero, see \cite{fukushima} for definitions of capacity, polar sets, etc. Note that for $u,v\in\F^D$ we can write
\begin{align}
\C(u,v)&=\frac 1 2 \int_D\int_D (u(x)-u(y))(v(x)-v(y))J_X(y-x)dx dy+\int_D
	u(x)v(x)\kappa_D(x) dx,\nonumber
\end{align}
where $\kappa_D(x)=\int_{D^c}J_X(y-x)dy$ is called the \emph{killing density} of
$X^D$. It is also the density of the killing measure from the Beurling-Deny representation of the Dirichlet form $(\C,\F^D)$, \cite[Section 3.2]{fukushima}. Furthermore, $(\C,\F^D)$ is a regular Dirichlet form on $L^2(D)$ with a special standard core $C_c^{\infty}(D)$.

By removing the killing part from $(\C,\F^D)$ we obtain a new bilinear form: for every $u,v\in C_c^\infty(D)$ let
$$\E(u,v)=\frac 1 2 \int_{D}\int_{D}(u(x)-u(y))(v(x)-v(y))J_X(y-x)dx\,dy.$$
By Fatou's lemma the symmetric form
$(\E,C_c^\infty(D))$ is closable in $L^2(D)$, i.e.~for every sequence $u_n\in C_c^\infty(D)$ such that $u_n\xrightarrow{L^2}0$, 
\[
\E(u_n-u_m,u_n-u_m)\xrightarrow{n,m\to\infty}0\,\,\Rightarrow\,\,\E(u_n,u_n)\xrightarrow{n\to\infty}0,
\]
so we take $\F$ to be the closure of $C_c^\infty(D)$ under the inner product $\E_1=\E+(\cdot,\cdot)_{L^2(D)}$. The closed symmetric form $(\E,\F)$ is Markovian since it operates on a normal contraction, i.e. for $u\in \F$ and $v\in L^2(D)$,
\[
|v(x)-v(y)|\leq |u(x)-u(y)|, \ |v(x)|\leq |u(x)|, \ \forall x,y \ \Rightarrow \ \E(v,v)\leq \E(u,u).
\]
Therefore, the form $(\E,\F)$ is a regular Dirichlet form. By Theorem \cite[Theorem 7.21]{fukushima} there exists a symmetric Hunt process $Y$
associated with $(\E,\F)$, taking values in $D$ with lifetime $\zeta$. Following \cite{cen}, we call $Y$ the \emph{censored (or resurrected) process} associated with $X$. Note that the censored process $Y$ can also be interpreted as the process obtained from the L\'evy process $X$ by restricting its jumping measure to $D\times D$. 

The following theorem is the analogue of \cite[Theorem 2.1]{cen} and provides two alternative constructions of the censored process, by using the Ikeda-Nagasawa-Watanabe piecing together procedure from \cite{INW} and by resurrection through a Feynman-Kac transform. The proof of the theorem is analogous to the proof in \cite{cen} and we refer the reader to that proof. 

\begin{Thm}
\label{thm:construction}
The following processes have the same distribution
\begin{enumerate}[(i)]
  \item The symmetric Hunt process $Y$ associated with the regular Dirichlet
  form $(\E,\F)$ on $L^2(D)$.
  \item The strong Markov process obtained from the symmetric Levy process $X^D$ in $D$ through the Ikeda-Nagasawa-Watanabe piecing together procedure.
  \item The process obtained from $X^D$ through the Feynman-Kac transform $e^{\int_0^t \kappa_D(X_s^D)ds}$.
\end{enumerate}
\end{Thm}

From the construction of the censored process $Y$ through the Ikeda-Nagasawa-Watanabe piecing together procedure it follows that the censored process $Y$ can be obtained from the symmetric L\'evy process $X$ by suppressing its jumps from $D$ to the complement $D^c$. Several useful properties of the censored process follow directly from Theorem \ref{thm:construction}.

\begin{Rem}\label{rem:construction}
\cite[Theorem 3.10.]{chungzao} and Theorem \ref{thm:construction}(iii) imply that, if $X$ has an absolutely continuous transition measure, then so does the corresponding censored process $Y$. Furthermore, the censored process $Y$ is irreducible.
\end{Rem}

From now on, assume that there exists a strictly positive non-increasing function $j$ satisfying \eqref{eq:j1} and \eqref{eq:j2}. Moreover, $j$ is a radial L\'evy density of an isotropic unimodal L\' evy process with characteristic exponent $\psi$, see \eqref{eq:j3}. Furthermore, assume that $\psi$ satisfies conditions $\Hj$ and $\Hd$. 

In order to investigate the boundary behaviour of the corresponding censored process we introduce a new process through its Dirichlet form. Let $({\mathcal E}^{\text{ref}},\F_\text{a}^{\text{ref}})$ be a Dirichlet form on $L^2(D)$ defined by 
\begin{align}
	\F_\text{a}^{\text{ref}}&=\left\{u\in L^2(D):\frac 1 2\int_{D}	\int_{D}(u(x)-u(y))^2J_X(y-x)dx\,dy<\infty\right\}\nonumber\\
	\E^{\text{ref}}(u,v)&=\frac 1 2\int_{D}\int_{D}(u(x)-u(y))(v(x)-v(y))J_X(y-x)dx\,dy,\quad u,v\in\F_\text{a}^{\text{ref}}\nonumber
\end{align}
and note that it is not necessarily regular on $D$. Using the trace theorem from \cite{wagner} for Besov-type spaces of generalized smoothness on $d$-sets, $d\leq n$, and the analogue of \cite[Theorem 2.2]{cen} we first show that $({\mathcal E}^{\text{ref}},\F_\text{a}^{\text{ref}})$ is the active reflected Dirichlet form associated with $(\E,\F)$ in the sense of Silverstein, see \cite[Theorem 6.2.13 and Section 6.3]{chen}. 

\begin{Def}\label{def:dset}
 A non-empty Borel set $D$ is called a $d$-set, $0<d\leq n$, if there exist positive constants $c_1$ and $c_2$ such that for all $x\in D$ and $r\in (0, 1]$,
 $$c_1 r^d\leq \mathcal H^d(D\cap B(x,r))\leq c_2 r^d.$$
 \end{Def}

\begin{Thm}{Trace theorem}\label{trace_theorem}\\
Let $D$ be a $n$-set in $\R^n$, $\C_1:=\C+(\cdot,\cdot)_{L^2}$ and $X$ a purely discontinuous symmetric L\'evy process such that \eqref{eq:j1}, \eqref{eq:j2}, \Hj and \Hd hold, then the normed space $(\F_\text{a}^{\text{ref}},\sqrt{{\mathcal E}^{\text{ref}}_1})$ is the restriction of space $(\F^{\R^n},\sqrt{\C_1})$ on $D$ in the following sense: there exist operators $R:\F^{\R^n}\to \F_\text{a}^{\text{ref}}$ and $E:\F_\text{a}^{\text{ref}}\to\F^{\R^n}$ such that
\begin{align*}
&Ru=u \text{ a.e. on }D\text{ and }{\mathcal E}^{\text{ref}}_1(Ru,Ru)\leq c_3\C_1(u,u),\quad\forall u\in\F^{\R^n}\\
&Eu=u \text{ a.e. on }D\text{ and }\C_1(Eu,Eu)\leq c_4{\mathcal E}^{\text{ref}}_1(u,u),\quad\forall u\in \F_\text{a}^{\text{ref}}
\end{align*}
for some constants $c_3, c_4>0$ and $REu=u$ a.e.~on $D$ for all $u\in \F_\text{a}^{\text{ref}}$. Operators $R$ and $E$ are called the continuous restriction and extension operator respectively.
\end{Thm}

Using Theorem \ref{trace_theorem}, proofs of the following results are analogous to ones in \cite{cen}, so we omit them here.

\begin{Thm}\label{thm:reflected}
Let $D$ be an open set in $\R^n$. The Dirichlet form $({\mathcal E}^{\text{ref}},\F_\text{a}^{\text{ref}})$ is the active reflected Dirichlet form associated with $(\E,\F)$, i.e.
\begin{align*}
&\F_a^{\text{ref}}=\{u\in L^2(D): u_k=((-k)\vee u)\wedge k\in \F_{\text{loc}} \mbox{ and } \sup_{k} \E^{\text{ref}}(u_k,u_k)<\infty\}\\
&\E^{\text{ref}}(u,u)=\lim\limits_{k\to\infty} \E^{\text{ref}}(u_k,u_k).
\end{align*}
Here $f\in \F_{\text{loc}}$ if for every relatively compact open set $D_0$ in $D$ there exists a function $f_0\in\F$ such that $f=f_0$ a.e. on $D_0$.
\end{Thm}

By \cite[Theorem 6.6.3]{chen} the active reflected Dirichlet form is the maximal Silverstein extension of the corresponding regular Dirichlet form. This means that the space $\F_b$ of bounded functions in $\F$ is an ideal in space $\F_{a,b}^{\text{ref}}$ of bounded functions in $\F_{a,b}^{\text{ref}}$, i.e.~$\F_b\subset \F_{a,b}^{\text{ref}}$ and $f g\in\F_b$ for every $f\in\F_b$, $g\in\F_{a,b}^{\text{ref}}$. Furthermore, by \cite[Theorem 6.6.5, Remark 6.6.7]{chen} a Dirichlet form $( \E^*,\F^*)$ is a Silverstein extension of a quasi-regular Dirichlet form $(\E,\F)$ on $L^2(D)$ if and only if there exists a symmetric Hunt process $Y^*$ associated with the Dirichlet form $( \E^*,\F^*)$ that extends $Y$ to some state space $D^*$ which contains $D$ as an $\E^*$-quasi-open subset of $D^*$ up to an $\E$-polar set. Therefore, there exists a compactification $D^*$ of $D$ such that the active reflected Dirichlet form $({\mathcal E}^{\text{ref}},\F_\text{a}^{\text{ref}})$ is regular on $L^2(D^*)$ and we call the corresponding process $Y^*$ the \emph{reflected process} associated with $Y$. The set $D^*\setminus D$ is Lebesgue negligible, but not necessarily of zero ${\mathcal E}^{\text{ref}}$-capacity. Since $\F$ is the $\E_1$-closure of $C_c^\infty(D)$, the process $Y^*$ killed upon leaving $D$ has the same distribution as $Y$. Using this correspondence between $Y$ and $Y^*$ we arrive to the analogue of \cite[Theorem 2.4]{cen}.

\begin{Thm}
 Let $D$ be an open set in $\R^n$ with finite Lebesgue measure and $\zeta$ the lifetime of process $Y$. The following statements are equivalent
 \begin{enumerate}[(i)]
  \item $\P_x (\zeta<\infty) > 0$ for some (and hence for all) $x \in D$;
  \item $\P_x (\zeta<\infty) =1$ for some (and hence for all) $x \in D$;
  \item $1\not\in\F$;
  \item $\F \neq \F_a^{\text{ref}}$.
 \end{enumerate}
 \end{Thm}

Let $D$ be an open $n$-set in $\R^n$. Since $C_c(\R^n)$ is the special standard core in $(\C,\F^{\R^n})$, by  Theorem \ref{trace_theorem} $C_c(\overline{D})\cap \F_a^{\text{ref}}$ is a core for $(\E^{\text{ref}},\F_a^{\text{ref}})$, and therefore $(\E^{\text{ref}},\F_a^{\text{ref}})$ is a regular Dirichlet form on $\overline{D}$. This means that we can take $D^*=\overline{D}$ and that $Y$ can be represented as the process $Y^*$ killed upon hitting the boundary $\partial D$.

\begin{Rem}\label{remark_polar}
Let $D$ be an open $n$-set. If $\F\subsetneq \F^{\text{ref}}_a$ then $Y$ is a proper subprocess of $Y^*$ and $\partial D$ is not $\E^{\text{ref}}$-polar. This implies that
\[
\P_x(Y_{\zeta-}\in\partial D,\zeta<\infty)>0,\quad\forall x\in D.
\]
Additionally, if $D$ has finite Lebesgue measure, $Y^*$ is recurrent and therefore $\zeta$ is finite almost surely and
\[
\P_x(Y_{\zeta-}\in\partial D,\zeta<\infty)=1,\quad\forall x\in D.
\]
Note that the aforementioned statements hold not only for q.e.~$x\in D$, but can also be extended for all $x\in D$. This is due to the fact that $Y$ has an absolutely continuous transition density, see Remark \ref{rem:construction}.
\end{Rem}
So we see that the question of boundary behaviour of the censored process $Y$ is related to $\E^{\text{ref}}$-polarity of the boundary $\partial D$. Since every compact set is of finite capacity, by \cite[Theorem 4.2.1]{fukushima} a set $A$ is $\E$-polar if and only if $\text{Cap}_Y(A)=0$. The same is true for $\C$-polar sets. Also, since $X$ and $Y$ have absolutely continuous transition densities (Remark \ref{rem:construction}), by \cite[Theorem 4.1.2]{fukushima} every $\E$-polar ($\C$-polar) set is polar for the process $Y$ ($X$). The proof of the following important characterizations of $\E^{\text{ref}}$-polar sets follows in the same way as in \cite[Theorem 2.5]{cen} and \cite[Corollary 2.6]{cen}.
\begin{Thm}\label{theorem_capacity}
Let $D$ be an open $n$-set in $\R^n$.
\begin{enumerate}[(i)]
\item A set $A\subset \overline{D}$ is $\E^{\text{ref}}$-polar if and only if it is polar for the process $X$.
\item A set $A\subset D$ is polar for the process $Y$ if and only if it is polar for the process $X$.
\item If $A\subset\partial D$ is polar for the process $X$ then
\[
\P_x(Y_{\zeta-}\in A)=0,\quad \forall x\in D.
\]
\end{enumerate}
\end{Thm}
The converse of Theorem \ref{theorem_capacity}(iii) is not true, for a counterexample see \cite[Remark 2.2]{cen}.

\begin{Cor}\label{cor:boundary}
Let $D$ be an open $n$-set in $\R^n$ and $\zeta$ lifetime of the censored process $Y$. Then the following statements are equivalent.
\begin{enumerate}[(i)]
\item $Y\neq Y^*$;
\item $\F\subsetneq \F_a^{\text{ref}}$;
\item $\partial D$ is not polar for process $X$;
\item $\P_x\left(\lim\limits_{t\uparrow\zeta}Y_t\in\partial D,\,\zeta<\infty\right)>0$ for every $x\in D$;
\item $\P_x\left(\lim\limits_{t\uparrow\zeta}Y_t\in\partial D,\,\zeta<\infty\right)>0$ for some $x\in D$.
\end{enumerate}
\end{Cor} 

We conclude this section with the proof of Theorem \ref{thm:boundary}, which partially answers the question of boundary behaviour of the censored L\'evy process $Y$ in terms of the scaling coefficients $\delta_1$, $\delta_2$, $\delta_3$ and the Hausdorff measure of the boundary $\partial D$.

\proofof {\bf Theorem \ref{thm:boundary}:}
Note that there exists a constant $c>1$ such that for every Borel set $A$ in $\R^n$
\[
c^{-1}\cp_{X^{(1)}}(A)\leq\cp_X(A)\leq  c\cp_{X^{(2)}}(A),
\]
where $X^{(i)}=(X^{(i)}_t)_{t\geq 0}$ is a symmetric $(2\delta_i)$-stable L\' evy process, $i=1,2$. Furthermore, recall that $X$ is recurrent if and only if $\int_{B(0,r)}\frac{1}{\psi(|\xi|)}d\xi=\infty$ for some $r>0$ and when $n=1$ all points are non-polar for $X$ if and only if $\int_1^\infty\frac{1}{\psi(x)}dx<\infty$. These conditions are satisfied when $\delta_3\geq\frac 1 2$ and $\delta_1>\frac 1 2$ respectively. The proof now follows from \cite[Theorem 2.7]{cen}.
\qed

\section{Generalized 3G theorem}\label{section:3G}

In this section we prove a stronger version of the 3G inequality for a purely discontinuous symmetric L\'evy process $X$ in bounded $\kappa$-fat open sets and the generalized 3G inequality. Suppose that \eqref{eq:j1}, \eqref{eq:j2} and \Hj hold and that $n\geq 2$. The 3G inequality will be essential in comparing the Green function of the censored process $Y$ killed outside of some ball $B\subset D$ with the Green function of the killed process $X^B$. Before we start with the proof of Theorem \ref{harmonic:3GThm}, we recall some basic potential-theoretical results that we use in the proof. 

 Let $B$ be a bounded open set in $\R^n$, $n\geq 2$. Using the two-sided estimates for the transition density $p^B_t(x,y)$ of $X^B$ the following lower and upper bound on the Green function $G_B(x,y)=\int_0^\infty p_t^B(x,y)dt$ for $X^B$ were obtained in \cite[Lemma 2.7, Lemma 2.8]{martin}.

\begin{Lem}\label{lowerupperG}
Let $R \in (0,1)$ and $B$ be a bounded open set such that diam$(B) \leq R$. The Green function $G_B(x,y)$ is finite and continuous on $B\times B\setminus d$ and
\begin{enumerate}[(i)]
\item there exists a constant $c_1=c_1(R,\psi,\gamma_1,\gamma_2)$ such that for all $x,y\in B$
$$G_B(x,y)\leq c_1\frac{\Phi(|x-y|)}{|x-y|^n},$$ 
\item for every $L>0$ there exists a constant $ c_2= c_2(L,R,\psi,\gamma_1,\gamma_2)> 0$ such that for all $x,y \in B$ with ${|x - y|\leq L(\delta_B(x)\wedge \delta_B(y)),}$
$$G_B(x,y)\geq c_2\frac{\Phi(|x-y|)}{|x-y|^n}.$$
\end{enumerate}
\end{Lem}
\begin{Def}
 Let $U$ be an open subset in $\R^n$. A Borel measurable function $u$ on $\R^n$ is harmonic on $U$ with respect to $X$ if
 \begin{equation*}
  u(x) = \mathbb E_x [u(X_{\tau_B} )],\quad x \in B,
 \end{equation*}
 for every bounded open set $B$ such that $\overline{B} \subset U$. A harmonic function $u$ is regular harmonic in $U$ if
 \begin{equation*}
  u(x) = \mathbb E_x [u(X_{\tau_U} ),\,\tau_U<\infty],\quad x \in U.
 \end{equation*}
 \end{Def}
The following scale invariant Harnack inequality and boundary Harnack principle for harmonic functions of purely discontinuous L\'evy processes were established in \cite[Theorem 2.2, Theorem 2.3(ii)]{martin}.

\begin{Thm}\label{HIX}
Let $L > 0$. There exists a constant $c_3 = c_3(L,\psi,\gamma_1,\gamma_2) > 1$ such that the following is true: If $x_1, x_2 \in \R^n$ and $r\in (0, 1)$ are such that $|x_1-x_2| < Lr$, then for every
non-negative function $h$ which is harmonic with respect to $X$ in $B(x_1,r)\cup B(x_2,r)$, we have
$$c_3^{-1} h(x_2)\leq h(x_1)\leq c_3 h(x_2).$$
\end{Thm}

\begin{Thm}\label{harmonic:BHPX}
There exists a positive constant $c_4=c_4(\psi,\gamma_1,\gamma_2)>0$ such that for every $x_0\in \R^n$, every open set $B \subset \R^n$, every $r \in (0, 1)$ and all non-negative functions $h, v$ in $\R^n$ which are regular harmonic in $B \cap B(x_0, r)$ with respect to $X$ and vanish a.e.~in $B^c \cap B(x_0, r)$, we have
$$\frac{h(x)}{v(x)}\leq c_4\frac{h(y)}{v(y)},\quad x,y\in B \cap B\left(x_0, \frac r 2\right).$$
\end{Thm}

The results in this chapter concern a special class of open sets called $\kappa$-fat sets.

\begin{Def}
An open set $D\subset \R^n$ is said to be \emph{$\kappa$-fat} if there exist some $R>0$ and $\kappa\in\left(0,\frac 1 2\right]$ such that for every $Q \in \partial D$ and $r\in(0,R)$ there exists a ball $B(A_r(Q),\kappa r) \subset D \cap \B(Q, r)$. The pair $(R,\kappa)$ is called the characteristics of the $\kappa$-fat open set $D$.
\end{Def}
Note that the ball of radius $r>0$ is a $\kappa$-fat open set with characteristics $\left(2r,\frac 1 2\right)$. Let $B$ be a bounded $\kappa$-fat open set with characteristics $(R,\kappa)$ and diam$(B)\leq r$, for some $r>0$. Fix $z_0\in B$ such that $\kappa R<\delta_B(z_0)\leq R$. By Lemma \ref{lowerupperG}(i) and \eqref{eq:phi} it follows that 
\[
G_B(x,z_0)\leq c_5\frac{\Phi(\delta_B(z_0))}{\delta_B(z_0)^n},\quad  x\in B\setminus B\left(z_0,\delta_B(z_0)/2\right)
\]
where $c_5=2^{n}c_1$ depends only on $r$, $\psi$, $\gamma_1$, $\gamma_2$ and $n$. Instead of working directly with the Green function $G_B$, we define a function $g_B$ on $B$ by
\begin{equation}\label{eq2}
g_B(x)=G_B(x,z_0)\wedge c_5\frac{\Phi(\delta_B(z_0))}{\delta_B(z_0)^n}
\end{equation}
and note that if $|x-z_0|>\frac{\delta_B(z_0)}2$ then $g_B(x)=G_B(x,z_0)$. Let $\varepsilon_1 = \frac {\kappa R}{24}$ and for $x,y\in B$ define $r(x,y)=\delta_B(x)\vee\delta_B(y)\vee|x-y|$ and
\begin{equation}\label{eq3}
\mathcal B(x,y)=\left\{\begin{array}{c l}
\left\{A\in B:\delta_B(A)>\frac{\kappa}2 r(x,y),\,|x-A|\vee|y-A|<5r(x,y)\right\}, & \text{if }r(x,y)<\varepsilon_1\\
\{z_0\}, & \text{if }r(x,y)\geq\varepsilon_1. 
 \end{array}\right.
\end{equation}

The proof of Theorem \ref{harmonic:3GThm} is divided into several parts. The first theorem follows the proof of \cite[Theorem 1.2]{twosided} and \cite[Theorem 2.4]{hansen}.

\begin{Thm}\label{harmonic:3GThm2}
There exists a constant $c_6 = c_6(r,a,\kappa,\psi,\gamma_1,\gamma_2)>1$ such that for every bounded $\kappa$-fat open set $B$ with characteristics $(R,\kappa)$ such that diam$(B)\leq r$ and $\frac{R}{\text{diam}(B)}\geq a$ and every $x,y\in B$ and $A \in \mathcal B(x, y)$,
\begin{equation}
c_6^{-1}\frac{g(x)g(y)\Phi(|x-y|)}{g(A)^2|x-y|^n}\leq G_{B}(x,y)\leq c_6\frac{g(x)g(y)\Phi(|x-y|)}{g(A)^2|x-y|^n},
\end{equation}
where $g=g_{B}$ and $\mathcal B(x,y)$ are defined by \eqref{eq2} and \eqref{eq3} respectively.
\end{Thm}

\proof
Let $r_0:=\frac 1 2(|x-y|\wedge \varepsilon_1)$. We only consider the case $\delta_B(x) \leq \delta_B(y) \leq \frac {\kappa r_0}2 $, that is case (g) in \cite{hansen}, which implies $r(x,y)=|x-y|$. The remaining cases follow analogously. 

Choose $Q_x,Q_y \in\partial B$ with $|Q_x - x| = \delta_B(x)$ and $|Q_y - y| = \delta_B(y)$ and let $x_1=A_{\frac{\kappa r_0}2}(Q_x)$ and $y_1=A_{\frac{\kappa r_0}2}(Q_y)$. This means that $x,x_1\in B\cap B(Q_x,\frac{\kappa r_0}2)$ and $y,y_1\in B\cap B(Q_y,\frac{\kappa r_0}2)$. Since 
$$|z_0-Q_x|\geq \delta_B(z_0) \geq \kappa R=24\varepsilon_1>r_0
\text{ and }
|y-Q_x|\geq |x-y|-\delta_B(x)\geq \left(2-\frac{\kappa}{2}\right)r_0> r_0$$
functions $G_B(\cdot,y)$ and $G_B(\cdot,z_0)$ are regular harmonic in $B\cap B(Q_x,\kappa r_0)$ and vanish outside $B$. Recall from \eqref{eq2} that
\begin{equation}\label{harmonic:g2}
\delta_B(z)<\frac{\delta_B(z_0)}2\quad\Rightarrow\quad g(z)=G_B(z,z_0).
\end{equation}
From the boundary Harnack principle, Theorem \ref{harmonic:BHPX} it follows that
$$c_4^{-1}\frac{G_B(x_1,y)}{g(x_1)}\leq \frac{G_B(x,y)}{g(x)}\leq c_4\frac{G_B(x_1,y)}{g(x_1)}.$$
On the other hand, since $|z_0-Q_y|> r_0$ and
$$|x_1-Q_y|\geq |x-Q_y|-|x_1-Q_x|-\delta_B(x)\geq \left(2-\frac{\kappa}{2}\right)r_0-\frac{\kappa r_0}2-\frac{\kappa r_0}2> r_0,$$
functions $G_B(x_1,\cdot)$ and $G_B(\cdot,z_0)$ are regular harmonic on $B\cap B(Q_y,\kappa r_0)$. Applying the boundary Harnack principle as before it follows that
$$c_4^{-1}\frac{G_B(x_1,y_1)}{g(y_1)}\leq \frac{G_B(x_1,y)}{g(y)}\leq c_4\frac{G_B(x_1,y_1)}{g(y_1)}.$$
By putting the two inequalities above together, we arrive to
$$c_4^{-2}\frac{G_B(x_1,y_1)}{g(x_1)g(y_1)}\leq \frac{G_B(x,y)}{g(x)g(y)}\leq c_4^2\frac{G_B(x_1,y_1)}{g(x_1)g(y_1)}.$$
Since $\delta_B(x_1)\wedge\delta_B(y_1)\geq\frac{\kappa^2 r_0}{2}$, $\varepsilon_1|x-y|\leq 2r_0\text{diam}(B)$ and 
\begin{equation}\label{harmonic:eq1}
 |x_1-y_1|\leq |x_1-x|+|x-y|+|y-y_1|<\kappa r_0+|x-y|+\kappa r_0\leq (1+\kappa)|x-y|
\end{equation}
it follows that $|x_1-y_1|\leq\frac{96(1+\kappa)}{\kappa^3a}(\delta_B(x_1)\wedge\delta_B(y_1))$, so by applying Lemma \ref{lowerupperG} we arrive to
$$\frac{c_3c_4^{-2}}{g(x_1)g(y_1)|x_1-y_1|^n\phi(|x_1-y_1|^{-2})}\leq \frac{G_B(x,y)}{g(x)g(y)}\leq \frac{c_1c_4^{2}}{g(x_1)g(y_1)|x_1-y_1|^n\phi(|x_1-y_1|^{-2})}.$$
Applying \eqref{eq:phi}, \eqref{harmonic:eq1} and $|x_1-y_1|\geq |x-y|-|x_1-x|-|y_1-y|\geq |x-y|-2\kappa r_0\geq\left(1-\kappa\right)|x-y|$ the previous inequality transforms to
$$\frac{c_3c_4^{-2}(1+\kappa)^{-n}(1-\kappa)^2}{g(x_1)g(y_1)|x-y|^n\phi(|x-y|^{-2})}\leq \frac{G_B(x,y)}{g(x)g(y)}\leq \frac{c_1c_4^{2}(1-\kappa)^{-n}(1+\kappa)^2}{g(x_1)g(y_1)|x-y|^n\phi(|x-y|^{-2})}.$$
Lastly, we have to show that for all $A\in\mathcal B(x,y)$
\begin{equation}
g(A)^2\asymp g(x_1)g(y_1).\label{show}
\end{equation}
Consider two cases, $r_0<\frac{\varepsilon_1}{2}$ and $r_0=\frac{\varepsilon_1}{2}$. If $r_0<\frac{\varepsilon_1}{2}$ then $r(x,y)=|x-y|<\varepsilon_1$, $r_0=\frac 1 2 r(x,y)$ and
$\delta_B(x_1)\wedge\delta_B(y_1)\geq \frac{\kappa^2 r(x,y)}{4}$. Since $G_B(\cdot,z_0)$ is harmonic on $B(x_1,\delta_B(x_1))\cup B(A,\delta_B(A))$ and 
\[
|x_1-A|\leq |x_1-x|+|x-A|\leq \kappa r_0+5r(x,y)\leq\frac{4}{\kappa^2}\left(\frac \kappa 2+5\right) (\delta_B(x_1)\wedge \delta_B(A))
\]
by \eqref{harmonic:g2} and Theorem \ref{HIX} it follows that $c_3^{-1} g(x_1)\leq g(A)\leq c_3 g(x_1)$. The analogous inequality follows for $y_1$ in place of $x_1$ and therefore $c_3^{-2} g(x_1)g(y_1)\leq g^2(A)\leq c_3^2 g(x_1)g(y_1)$. 

On the other hand, if $r_0=\frac{\varepsilon_1}{2}$ then $r(x,y)=|x-y|\geq \varepsilon_1$, so by \eqref{eq3} and \eqref{eq2} it follows that 
\begin{equation}\label{showa}
g(A)=g(z_0)=c_5 \frac{\Phi(\delta_B(z_0))}{\delta_B(z_0)^n}.
\end{equation}
Let $v\in\{x_1,y_1\}$ and $z\in B$ such that $|z-z_0|=\frac{\delta_B(z_0)}2=\delta_B(z)$. Since $\delta_B(v)\geq \frac{\kappa^2\varepsilon_1}{4}$ it follows that $|v-z|\leq \text{diam}(B)\leq \frac{96}{\kappa^3 a}(\delta_B(v)\wedge\delta_B(z))$, so by applying Theorem \ref{HIX} we get
$$c_3^{-1} G_B(z,z_0)\leq g(v)\leq c_3^{-1} G_B(z,z_0).$$
Therefore, by Lemma \ref{lowerupperG} and \eqref{eq:phi} it follows that
\[
\tilde c^{-1} \frac{\Phi(\delta_B(z_0))}{\delta_B(z_0)^n}\leq g(v)\leq \tilde c\frac{\Phi(\delta_B(z_0))}{\delta_B(z_0)^n}
\]
for some $\tilde c=\tilde c(r,a,\kappa,\psi,\gamma_1,\gamma_2)>1$, which together with \eqref{showa} implies \eqref{show}.
\qed

We will also need the following result from \cite[Lemma 2.7]{minimal}.

\begin{Lem}\label{carleson}{\bf Carleson's estimate}\\
Let $r>0$ and $\kappa\in\left(0,\frac 1 2\right]$. There exists a constant $c_7 = c_7(r,\kappa,\psi,\gamma_1,\gamma_2)>0$ such that for every bounded open $\kappa$-fat set $B$ with characteristics $(R,\kappa)$ and diam$(B)\leq r$, $z\in\partial B$, $r_0\in(0,\frac{\kappa R}{4})$ and $y\in B\setminus \overline{B(z,3r_0)}$
$$G_{B}(x,y)\leq c_7  G_{B}(A_{r_0}(z),y),\quad x\in B\cap B(z,r_0).$$

 \end{Lem}

Applying the Carleson's estimate, the Harnack inequality and Lemma \ref{lowerupperG} the proofs of the following lemmas follow entirely as in \cite[Lemma 3.8-3.11]{generalized3G}. Let $B$ be a bounded $\kappa$-fat open set with diam$(B)\leq r$ and characteristics $(R,\kappa)$ such that $\frac{R}{\text{diam}(B)}\geq a$. As in the proof of Theorem \ref{harmonic:3GThm2}, let $Q_x\in\partial B$ be such that $|x-Q_x|=\delta_B(x)$, $x\in B$.

\begin{Lem} There exists a constant $c_8=c_8(r,a,\kappa,\psi,\gamma_1,\gamma_2)>0$ such that for every $x,y\in B$ with $r(x,y)<\varepsilon_1$,
\begin{equation}
g(z)<c_8 g(A_{r(x,y)}(Q_x)),\quad z\in B\cap B(Q_x,r(x,y)).\label{g1}
\end{equation}
\end{Lem}

\begin{Lem} There exists a constant $c_9=c_9(r,a,\kappa,\psi,\gamma_1,\gamma_2)>0$ such that for every $x,y\in B$ 
\begin{equation}
g(x)\vee g(y)<c_9g(A),\quad A\in \mathcal B(x,y).\label{g2}
\end{equation}

\end{Lem}

\begin{Lem}If $x,y,z\in B$ satisfy $r(x,z)\leq r(x,y)$, then there exists a constant $c_{10}=c_{10}(r,a,\kappa,\psi,\gamma_1,\gamma_2)>0$ such that
\begin{equation}
g(A_{x,y})<c_{10} g(A_{y,z}),\quad \text{for every }(A_{x,y},A_{y,z})\in \mathcal B(x,y)\times \mathcal B(y,z).\label{g3}
\end{equation}
\end{Lem}

\begin{Lem} There exists a constant $c_{11}=c_{11}(r,a,\kappa,\psi,\gamma_1,\gamma_2)>1$ such that for every $x,y,z,w\in B$ and $(A_{x,y},A_{y,z},A_{x,z})\in \mathcal B(x,y)\times \mathcal B(y,z)\times B(x,z)$,
\begin{equation}
g(A_{x,z})^2<c_{11}\left(g(A_{x,y})^2+g(A_{y,z})^2\right)\label{g4}
\end{equation}
\end{Lem}
\proofof {\bf Theorem \ref{harmonic:3GThm}:} Applying Theorem \ref{harmonic:3GThm2} we get
$$\frac{G_B(x,y)G_B(y,z)}{G_B(x,z)}\leq c_6^3 \frac{g(y)^2g(A_{xz})^2}{g(A_{xy})^2g(A_{yz})^2}\frac{\Phi(|x-y|)\Phi(|y-z|)}{\Phi(|x-z|)}\frac{|x-z|^n}{|x-y|^n|y-z|^n}.$$
By \eqref{g4} and \eqref{g2}, 
\begin{align*}
\frac{g(y)^2g(A_{xz})^2}{g(A_{xy})^2g(A_{yz})^2}&\leq c_{11} \left(\frac{g(y)^2}{g(A_{xy})^2}+\frac{g(y)^2}{g(A_{yz})^2}\right)\leq 2c_{11} c_9^2,
\end{align*}
which proves the 3G inequality \eqref{harmonic:3G} with $c_1=2c_6^3c_{11} c_9^2$ depending only on $r$, $a$, $\kappa$, $\psi$, $\gamma_1$ and $\gamma_2$.
\qed

Next we show the generalized 3G inequality, following the approach in \cite{generalized3G}

\begin{Lem}\label{lem4} There exist positive constants $c_{12}=c_{12}(R,\kappa,\psi,\gamma_1,\gamma_2)$ and $\beta=\beta(R,\kappa,\psi,\gamma_1,\gamma_2)<2\delta_2$ such that for every bounded $\kappa$-fat open set $D$, $Q\in\partial D$, $r\in(0,R)$ and non-negative function $u$ which is harmonic with respect to $X$ in $D\cap B(Q,r)$ we have
\[
u(A_r(Q))\leq c_{12}\left(\frac r s\right)^\beta u(A_s(Q)),\quad s\in(0,r).
\]
\end{Lem}
\proof
Here we follow the proof of \cite[Lemma 5.2]{BHPforSBM}. For $k\in\mathbb N_0$ and $r\in(0,R)$ let $\eta_k:=\left(\frac \kappa 2\right)^k r$, $A_k:=A_{\eta_k}(Q)$ and $B_k:=B(A_k,\eta_{k+1})$. Since $B(A_k,2\eta_{k+1})\subset D$ it follows that $B(Q,\eta_{k+1})\cap B_k=\emptyset$ so the balls $B_k$ are disjoint. By the harmonicity of $u$ and Theorem \ref{HIX} we have
\begin{equation}\label{lemma:show}
u(A_k)\geq =\sum_{i=0}^{k-1}\int_{B_i} u(y)K_{B_k}(A_k,y)dy\geq c_3^{-1}\sum_{i=0}^{k-1}u(A_i)\int_{B_i}K_{B_k}(A_k,y)dy.
\end{equation}
By \cite[Lemma 2.9]{martin} and \eqref{eq:j4} there exist constants $\tilde c,\tilde c_1>0$ such that for $i\in\{0,...,k-1\}$ and $y\in B_i$
\[
K_{B_k}(A_k,y)\geq \tilde c\frac{j(|y-A_k|)}{\phi(\eta_{k+1}^{-2})}\geq \tilde c_1\frac{\phi(\eta_i^{-2})}{\eta_i^{n}\phi(\eta_{k+1}^{-2})}
\]
which together with \eqref{lemma:show} implies that $u(A_k)\geq \tilde c_2 \sum_{i=0}^{k-1}u(A_i)\frac{\phi(\eta_i^{-2})}{\phi(\eta_{k}^{-2})}$, for some $\tilde c_2>0$. Iterating this inequality we get that $u(A_k)\phi(\eta_{k}^{-2})\geq \tilde c_2(1+\tilde c_2)^{k-1}u(A_0)\phi(\eta_0^{-2})$ for every $k\in\N$.
This inequality together with \Hj and \eqref{eq:phipsi} implies that  
\begin{align}
u(A_r(Q))&\leq \frac {a_2\gamma_2}{\tilde c_2(1+\tilde c_2)^{k-1}}\left(\frac{2}{\kappa}\right)^{2k\delta_2}u(A_{k})=\frac {a_2(1+\tilde c_2)}{\tilde c_2}\left(\frac{2}{\kappa}\right)^{k\beta}u(A_{k}),\label{eq4}
\end{align}
where $\beta=2\delta_2-\frac{\log{(1+\tilde c_2)}}{\log{\left(\frac 2 \kappa\right)}}$. Let $0<s<r$ and $k\in\N$ such that $\eta_k\leq s<\eta_{k-1}$. Since $B(A_s(Q),\kappa\eta_k)\cup B(A_k,\kappa\eta_k)\subset B$ and $|A_k-A_s(Q)|\leq 2\eta_{k-1}=\frac 4\kappa\eta_k$, by Theorem \ref{HIX} the inequality \eqref{eq4} transforms to
\begin{align*}
u(A_r(Q))\leq \frac {c_3a_2(1+\tilde c_2)}{\tilde c_2}\left(\frac{2}{\kappa}\right)^{\beta}\left(\frac{r}{s}\right)^{\beta}u(A_{s}(Q))=\tilde c_3 \left(\frac{r}{s}\right)^{\beta}u(A_{s}(Q)). 
\end{align*}
\qed

We will use the following remark several times in the proofs of the following lemmas.

\begin{Rem}\label{rem}
 Let $x,y\in B$ such that $r(x,y)<\varepsilon_1$. As in the proof of Theorem \ref{harmonic:3GThm2} let $Q_x,Q_y\in\partial B$ be such that $|x-Q_x|=\delta_B(x)$ and $|y-Q_y|=\delta_B(y)$. Note that $A_{r(x,y)}(Q_x),A_{r(x,y)}(Q_y)\in \mathcal B(x,y)$, since for $A:=A_{r(x,y)}(Q_x)$ it follows that $\delta_B(A)\geq \kappa r(x,y)$,
$|x-A|\leq \delta_D(x)+|Q_x-A|\leq 2r(x,y)$ and $|y-A|\leq |y-x|+|x-A|\leq 3r(x,y).$ By Theorem \ref{HIX}, $g(A_1)\asymp g(A_2)$ for all $A_1,A_2\in\mathcal B(x,y)$ and therefore it follows that
\[
 g(A_{x,y})\asymp g(A_{r(x,y)}(Q_x))\asymp g(A_{r(x,y)}(Q_y)),
\]
for all $A_{x,y}\in\mathcal B(x,y)$.
\end{Rem}

\begin{Lem}\label{lem5}
There exists a constant $c_{13}=c_{13}(r,\kappa, R,\psi,\gamma_1,\gamma_2)>0$ such that for every $\kappa$-fat open set $B$ with characteristics $(\kappa,R)$ and diam$(B)\leq r$ and every $x,y\in B$ with $r(x,y)<\varepsilon_1$
\[
 g(A_{x,y})\geq c_{13} r(x,y)^\beta,
\]
for all $A_{x,y}\in\mathcal B(x,y)$.
\end{Lem}
\proof
By Remark \ref{rem} it is enough to prove the inequality for $A=A_{r(x,y)}(Q_x)$. Since $\delta_B(z_0)>\kappa R=24\varepsilon_1$, function $g$ is harmonic in $B\cap B(Q_x,\varepsilon_1)$ so by Lemma \ref{lem4} it follows that
\begin{equation}\label{eq5}
 g(A)\geq c_{12}^{-1}\left(\frac{r(x,y)}{\varepsilon_1}\right)^\beta g(A_{\varepsilon_1}(Q_x)).
\end{equation}
Since $\frac{\delta_B(z_0)}{2}\leq |A_{\varepsilon_1}(Q_x)-z_0|\leq \frac{2r}{24\varepsilon_1}\delta_{B}(z_0)$, by Lemma \ref{lowerupperG} it follows that $g(A_{\varepsilon_1}(Q_x))\geq c_2\frac{\Phi(|A_{\varepsilon_1}(Q_x)-z_0|)}{|A_{\varepsilon_1}(Q_x)-z_0|^n}\geq c_2\frac{\Phi\left(\frac{\kappa R}{2}\right)}{(2r)^n}>0$. This inequality together with \eqref{eq5} completes the proof.
\qed

\begin{Lem}\label{lem7}
There exists a constant $c_{14}=c_{14}(r,\kappa,R,\psi,\gamma_1,\gamma_2)>0$ such that for every $\kappa$-fat open set $B$ with characteristics $(R,\kappa)$ and all $x,y,z\in B$ and $(A_{x,y},A_{y,z})\in\mathcal B(x,y)\times \mathcal B(y,z)$
 \[
  \frac{g(A_{y,z})}{g(A_{x,y})}\leq c_{14}\left[\left(\frac{r(y,z)}{r(x,y)}\right)^\beta\vee 1\right].
 \]
\end{Lem}
 \proof
 The proof follows the proof of \cite[Lemma 3.13]{generalized3G}, by applying Lemma \ref{carleson}, Lemma \ref{lem4}, Lemma \ref{lem5} and Remark \ref{rem}.  
\qed

Finally, the generalized 3G inequality now follows by adapting the arguments of \cite[Theorem 1.1]{generalized3G}. Let
\[
H(x,y,z,w)=\frac{\Phi(|x-y|)\Phi(|z-w|)}{\Phi(|x-w|)}\frac{|x-w|^n}{|x-y|^n|z-w|^n}.
\]

\begin{Thm}{(Generalized 3G theorem)}\label{Generalized3GThm}\\
Let $r>0$, $R>0$ and $\kappa\in\left(0,\frac 1 2 \right]$. There exist constants $\beta=\beta(R,\kappa,\psi,\gamma_1,\gamma_2)\leq 2\delta_2$ and $c_{15}=c_{15}(r,R,\kappa,\psi,\gamma_1,\gamma_2)>0$ such that 
\begin{equation*}
\frac{G_B(x,y)G_B(z,w)}{G_B(x,w)}\leq c_{15}\left(\frac{|x-w|\wedge |y-z|}{|x-y|}\vee 1\right)^\beta\left(\frac{|x-w|\wedge |y-z|}{|z-w|}\vee 1\right)^\beta H(x,y,z,w) 
\end{equation*}
for every bounded $\kappa$-fat open set $B$ with characteristics $(R,\kappa)$ and diam$(B)\leq r$.
\end{Thm}

\section{Harnack inequality for censored L\' evy processes}\label{section:harnack}

Let $D$ be an open set in $\R^n$, $n\geq 2$ and suppose that \eqref{eq:j1}, \eqref{eq:j2}, \eqref{eq:j3} hold and that the characteristic exponent $\psi$ satisfies the scaling condition \Hj. As before, let $X$ and $Y$ be the symmetric pure jump L\'evy process with characteristic exponent $\Psi_X$ and the corresponding censored process on $D$, respectively. Recall that, by Theorem \ref{thm:construction}, $Y$ can be obtained through the Feynman-Kac transform by resurrecting $X^D$ at the rate $\kappa_D$. Let $B$ be a bounded Borel set in $D$ and $p_t^{B}$ the transition density of $X^B$. From \cite[Lemma 3.5]{gauge} it follows that the Green function $G_B^Y$ for $Y^B$ is equal to
\begin{equation}\label{eq:gauge}
G_B^Y(x,y)=u(x,y)G_B(x,y),\, x,y\in B,
\end{equation}
where $G_B$ is the Green function for $X^B$ and 
$$u(x,y)=\EE_x^y\left[e^{A(\tau_B)}\right]=\EE_x^y\left[e^{\int_0^{\tau_B}\kappa_D(X_s)ds}\right]$$
is called the conditional gauge function. Here we denote by $\P_x^y$ and $\EE_x^y$ the probability and expectation for the $G_B(x,y)$-conditioned process starting from $x\in B$ respectively, i.e. the process with transition density
\[
p^{y}_t(x,z)=\frac{G_B(x,z)}{G_B(x,y)}p_t^{B}(x,z),\quad t>0,\,x,z\in B.
\]
In order to obtain two-sided estimates for the Green function $G_B^Y$ on a ball $B\subset D$, it suffices to show that the conditional gauge function $u$ is bounded, i.e. that the Green functions $G_B$ and $G_B^Y$ are comparable. The following result is the analogue of \cite[Lemma 3.1]{cen}. Due to the lack of exact scaling, we use a more general version of the 3G inequality showed in the previous section, Theorem \ref{harmonic:3GThm}.

\begin{Lem}\label{harmonic:Harnack_lemma}
 There is a constant $r_1=r_1(\psi,\gamma_1,\gamma_2)\in(0, \frac 1 3)$, independent of $D$, such that for every $r\in(0,1)$ and every ball $B=B(x,r_1r)\subset B(x,r)\subset D$,
 $$\int_B \frac{G_B(v,y)G_B(y,w)}{G_B(v,w)}\kappa_D(y)dy\leq\frac 1 2,\quad\forall v,w\in B.$$
\end{Lem}
\proof
Let $r_1\leq \frac 1 3$ and $r\in(0,1)$. Since $rr_1<\frac 1 3$ by Theorem \ref{harmonic:3GThm}  
\begin{align*}
 \frac{G_B(v,y)G_B(y,w)}{G_B(v,w)}&\leq c_1  \frac{\phi(|v-w|^{-2})}{\phi(|v-y|^{-2})\phi(|y-w|^{-2})}\frac{|v-w|^n}{|v-y|^n|y-w|^n},\quad\forall v,y,w\in B.
\end{align*}
First we will show that there exists a constant $\tilde c=\tilde c(n,\phi)>0$ such that
$$\phi(|v-w|^{-2})|v-w|^n\leq \tilde c\left(\phi(|v-y|^{-2})|v-y|^{n}+\phi(|y-w|^{-2})|y-w|^{n}\right).$$
From \eqref{eq:phi} it follows that $ \phi(s^{-2})s^n\leq \frac{r^{2}}{s^{2}}\phi(r^{-2})s^{n}\leq  \phi(r^{-2})r^n$, for all $s<r\leq1.$ Without loss of generality, we can assume that $|v-y|\leq |y-w|$, so $|v-w|\leq 2|y-w|$ and
\begin{align*}
 \phi(|v-w|^{-2})|v-w|^n\leq 2^{n}\left(\phi(|v-y|^{-2})|v-y|^n+\phi(|y-w|^{-2})|y-w|^n\right). 
\end{align*}
Therefore for every $v,w\in B$
\begin{align*}
\int_B\frac{G_B(v,y)G_B(y,w)}{G_B(v,w)}dy&\leq c_12^{n}\left( \int\limits_{B(v,2rr_1)}\frac 1 {\phi(|v-y|^{-2})|v-y|^{n}}dy + \int\limits_{B(w,2rr_1)}\frac 1{\phi(|y-w|^{-2})|y-w|^{n}}dy\right)\\
& \leq \tilde c_1 \int_0^{2rr_1}\frac 1 {\phi(s^{-2})s^{n}} s^{n-1}ds\overset{\Hj}{\leq} \frac{\tilde c_1 }{2a_1\gamma_2\delta_1  }\phi((2r_1r)^{-2})^{-1},
\end{align*}
for some $\tilde c_1=\tilde c_1(\psi)>0$. Furthermore, for every $y\in B=B(x,rr_1)\subset D$ it follows that $B(y,r(1-r_1))\subset D$, which by \cite[Lemma 2.2]{martin} implies that 
\begin{align*}
\kappa_D(y)&\leq\tilde c_2\int_{r(1-r_1)}^\infty s^{n-1}j(s)ds\leq\tilde c_3 \phi(r^{-2}(1-r_1)^{-2}),
\end{align*}
for all $r >0$ and constants $\tilde c_2, \tilde c_3>0$ depending only on $\psi$, $\gamma_1$ and $\gamma_2$. Finally, for $r_1$ small enough we get the following,
\begin{align*}
\int_B \frac{G_B(v,y)G_B(y,w)}{G_B(v,w)}\kappa_D(y)dy&\leq \frac{\tilde c_1\tilde c_3}{2a_1\gamma_2\delta_1} \frac{\phi(r^{-2}(1-r_1)^{-2})}{\phi((2r_1r)^{-2})}\overset{\Hj}{\leq}  \frac{\tilde c_1\tilde c_3}{2a_1^2\gamma_2^2\delta_1}\left(\frac{2r_1}{1-r_1}\right)^{2\delta_1}\leq \frac 1 2. 
\end{align*}
\qed

By the previous lemma it follow that for every $r<1$, every ball $B=B(x,r r_1)\subset D$ and $v,w\in B$ 
\begin{align*}
\EE_v^w[A(\tau_B)]&=\int_B\kappa_D(y)\frac{G_B(v,y)G_B(y,w)}{G_B(v,w)}dy\leq \frac 1 2,
\end{align*}
so by Khasminskii's lemma, \cite[Lemma 3.7]{chungzao}
\begin{equation}
1\leq u(v,w)=\EE_v^w[e^{A(\tau_B)}]\leq\frac 1{1-\frac 1 2}=2.\label{harmonic:khasminskii}
\end{equation}

Let $B$ be a ball such that $\overline{B}\subset D$. By calculations in \cite[p.318]{martin}, process $X$ satisfies the conditions of \cite[Theorem 1]{sztonyk} which implies that for all $y\in B$
 $$\P_y (X_{\tau_{B}} \in \partial B) = \P_y\left(X_{\tau_{B}-}= X_{\tau_{B}}\right)=0.$$
Note that by Theorem \ref{thm:construction}(iii) it follows that $\P_y\left(Y_{\tau^Y_{B}-}= Y_{\tau^Y_{B}}\right)=0$, for all $y\in B\subset\overline{B}\subset D.$ Using the L\' evy system formula and \eqref{eq:gauge} we arrive to the formula for the joint distribution of $(Y_{\tau^Y_B-} , Y_{\tau^Y_B} )$ restricted to the event $\{\tau_B^Y<\infty\}$, 
\begin{equation}
\EE_x[f(Y_{\tau^Y_B-})g(Y_{\tau^Y_B})]=\int_{D\setminus\overline{B}}\int_B f(y)g(z)G_B(x,y)u(x,y)J_X(z-y)dy dz,\label{harmonic:IWformula}
\end{equation}
for all non-negative Borel measurable functions $f$ and $g$ on $D$ and open Borel sets $B\subset\overline{B}\subset D$. Recall that the Poisson kernel $K_B$, i.e. the density function of the $\P_x$-distribution of $X_{\tau_B}$, is of the form
\begin{equation*}
K_B(x,z)=\int_B G_B(x,y)J_X(z-y)dy,\quad x\in B,\,z\in B^c.
\end{equation*} 

Using \eqref{harmonic:khasminskii} and \eqref{harmonic:IWformula} we are able to prove the scale invariant Harnack inequality for harmonic functions with respect to the censored process $Y$.

 \begin{Thm}
  \label{harnack inequality}
  For any $L > 0$, there exists a constant $c_{1} = c_{1}(\psi,\gamma_1,\gamma_2, L) > 1$ such that the following is true: If $x_1,x_2\in D$ and $r\in(0,1)$ are such that $B(x_1,r)\cup B(x_2,r)\subset D$ and $|x_1-x_2 | < Lr$, then for every
non-negative function $h$ which is harmonic with respect to $Y$ on $B(x_1,r)\cup B(x_2,r)$, we have
  $$c_1^{-1} h(x_1)\leq h(x_2) \leq c_1 h(x_1).$$
 \end{Thm}

\proof
Let $r_1\in(0,\frac 1 3)$ be the constant from Lemma \ref{harmonic:Harnack_lemma} and $B_i=B(x_i,r_1r),$ $i=1,2$. Since $\overline{B_1}\subset D$ it follows that for $y\in B_1$
\begin{align*}
h(y)&=\mathbb E_y\left[h(Y_{\tau_{B_1}^Y})\right]\overset{\eqref{harmonic:IWformula}}{=}\int_{D\setminus\overline{B}_1}\int_{B_1}h(w)G_{B_1}(y,v)u(y,v)J_X(w-v)dv\,dw=\EE_y\left[h(X_{\tau_{B_1}})u(y,X_{\tau_{B_1}-})\right].
\end{align*}
Here we implicitly assume $h=0$ on $D^c$. Define $w(y):=\EE_y\left[h(X_{\tau_{B_1}})\right]$, $y\in B_1$, and note that $w$ is harmonic in $B_1$ with respect to $X$. From \eqref{harmonic:khasminskii} it follows that 
 \begin{equation}\label{harmonic:w}
 w(y)\leq h(y)\leq 2w(y),\quad\forall y \in B_1 
 \end{equation}
 and analogously
 \begin{equation}\label{harmonic:w2}
 \EE_y\left[h(X_{\tau_{B_2}})\right]\leq h(y)\leq 2\EE_y\left[h(X_{\tau_{B_2}})\right],\quad\forall y \in B_2. 
 \end{equation} 
 By \cite[Proposition 2.3]{martin} there exists a constant $\tilde c_1=\tilde c_1(\psi,\gamma_1,\gamma_2)>0$ such that for any $y\in B(x_1,\frac{rr_1}{2})$
 \begin{align}\label{harmonic:w3}
  w(y)=\int_{D\setminus \overline{B}_1} h(z)K_{B_1}(y,z)dz\geq \tilde c_1\int_{D\setminus \overline{B}_1} h(z)K_{B_1}(x_1,z)dz=\tilde c_1 w(x_1)\geq \frac{\tilde c_1}2 h(x_1).
 \end{align}
First we consider the case when $r\leq |x_1-x_2|<Lr$. It follows that $B_2\cap B(x_1,r_1r/2)=\emptyset$ and therefore by \eqref{harmonic:w}, \eqref{harmonic:w2} and \eqref{harmonic:w3} 
 \begin{align}
  h(x_2)&\geq\EE_{x_2}\left[w(X_{\tau_{B_2}});X_{\tau_{B_2}}\in B(x_1,r_1r/2)\right]\geq\frac{\tilde c_1}2 h(x_1)\int_{B(x_1,r_1r/2)} K_{B_2}(x_2,z)dz.\label{harmonic:HIproof1}
 \end{align}
By \cite[Lemma 2.6]{martin} there exists a constant $\tilde c_2=\tilde c_2(\psi,\gamma_1,\gamma_2)>0$ such that for all $z\in\overline{B}_2^c$
\begin{equation}\label{harmonic:poisson}
K_{B_2}(x_2,z)\geq \tilde c_2\frac{ j(|z-x_2|)}{\phi((r_1r)^{-2})}.
\end{equation} 
Also, for $z\in B(x_1,r_1r/2)$ we have $|z-x_2|\leq r(r_1/2+L)<r_1/2+L$, so by \eqref{eq:j4} there exists a constant $\tilde c_3=\tilde c_3(\psi, \gamma_1,\gamma_2, L)>0$ such that
\begin{equation}\label{harmonic:j}
 j(|z-x_2|)\geq j(r(r_1/2+L))\geq \tilde c_3 \frac{\phi(r^{-2}(r_1/2+L)^{-2})}{r^n(r_1/2+L)^n}.
\end{equation}
Combining \eqref{harmonic:HIproof1}, \eqref{harmonic:poisson} and \eqref{harmonic:j} we get
 \begin{align*}
  h(x_2)&\geq\frac{\tilde c_1\tilde c_2\tilde c_3}{2}\frac{|B(x_1,\frac{rr_1}{2})|}{r^n (r_1/2+L)^n}\frac{\phi(r^{-2}(r_1/2+L)^{-2})}{\phi((r_1r)^{-2})}h(x_1)\\
&\overset{\eqref{eq:phi}}{\geq}\frac{\tilde c_1\tilde c_2\tilde c_3}{2}\frac{|B(x_1,\frac{rr_1}{2})|}{r^n (r_1/2+L)^n}\left(1\wedge \left(\frac{r_1}{\frac{r_1}{2}+L}\right)^2\right)h(x_1)= c_{1}(\psi,\gamma_1,\gamma_2,L) h(x_1).
\end{align*}
On the other hand, if $|x_1-x_2|<r$ take $r'=|x_1-x_2|$ and $L'=1$. Since $r'\leq |x_1-x_2|<L'r'$ the proof follows in the same way as in the previous case.  
\qed

\begin{Rem}
If for a Lipschitz domain $B\subset \overline{B}\subset D$ 
\[
 \inf\limits_{y\in B}\int_{D\setminus B} j(|z-y|)dz\geq c
\]
for some constant $c>0$, then by \eqref{harmonic:IWformula} it follows that
\[
 1=\int_{D\setminus \overline{B}}\int_B G_B^Y(x,y)J_X(z-y)dy\, dz\geq c \int_B G_B^Y(x,y)dy
\]
and therefore
\begin{equation}\label{harmonic:rem}
 \EE_x[\tau_B^Y]=\int_{B}G_B^Y(x,y)dy<\infty,\quad \forall x\in B.
\end{equation}
Furthermore, \eqref{harmonic:rem} holds for all $x \in D$ and implies that
\[
 \P_x(\tau_B^Y<\infty)=1, \quad\text{ for all }x\in D.
\]
\end{Rem}

\vspace{1cm}

\bibliographystyle{plain}
\bibliography{censored_bib}

\end{document}